\documentclass{article}

\usepackage[letterpaper,margin=1in]{geometry}
\usepackage[latin1]{inputenc}
\usepackage[T1]{fontenc}
\usepackage[english]{babel}
\usepackage{amsmath,amsthm,amsfonts}
\usepackage{url}

\makeatletter
\renewcommand\section{\@startsection {section}{1}{\z@}%
                                   {-3.5ex \@plus -1ex \@minus -.2ex}%
                                   {2.3ex \@plus.2ex}%
                                   {\normalfont\Large\bfseries}}
\renewcommand\subsection{\@startsection{subsection}{2}{\z@}%
                                     {-3.25ex\@plus -1ex \@minus -.2ex}%
                                     {1.5ex \@plus .2ex}%
                                     {\normalfont\large\itshape\bfseries}}
\makeatother

\newtheorem{theorem}{\large Theorem} 
\newtheorem{proposition}[theorem]   {\large Proposition}

\newtheorem{corollary}[theorem]{\large  Corollary}

\newtheorem{remark}[theorem]{\large Remark}

\font\tenmath=msbm10

\newfam\mathfam \textfont\mathfam=\tenmath

\def \\ { \cr }

\newcommand{\EE}{{\mathbb E}}

\newcommand{\NN}{{\mathbb N}}
\newcommand{\PP}{{\mathbb P}}

\newcommand{\A}{{\cal A}}

\newcommand{\aD}{{\cal D}}
\newcommand{\eL}{{\cal L}}

\newcommand{\X}{{\cal X}}
\newcommand{\Y}{{\cal Y}}
\newcommand{\Z}{{\cal U}}
\newcommand{\hX}{{{\widehat{\cal X}}}}
\newcommand{\hZ}{{{\widehat{Z}}}}

\newcommand{\hI}{{{\widehat{I}}}}
\newcommand{\tM}{{{\widehat{M}}}}

\newcommand{\tth}{{{\widehat{h}}}}
\newcommand{\hnu}{{{\widehat{\nu}}}}

\newcommand{\aI}{{\cal I}}

\newcommand{\fe}{\mathfrak{e}}
\newcommand{\bp}{{\mathfrak{b}}{\mathfrak{p}}}

\newcommand{\bz}{{\mathfrak{z}}}

\newcommand{\hxi}{{\widehat{\xi}}}

\newcommand\norm[1]{{\Vert{#1}\Vert}}
\newcommand\zeros{\ensuremath{{\mathbf0}}}
\newcommand\ones{\ensuremath{{\mathbf1}}}
\newcommand\ind{{\bf1}}

\begin{document}

\pagestyle{myheadings}
\markright{simple finite delayed multi-type branching process}

\title{A simple finite delayed multi-type branching process for 
infectious disease modeling}
\author{Andrew Hart%
\thanks{E-mail: \texttt{ahart@dim.uchile.cl}}%
\qquad Servet Mart{\'\i}nez%
\thanks{E-mail: \texttt{smartine@dim.uchile.cl}}
\\ \\
\parbox{5in}{\rmfamily\upshape\normalsize
Center for Mathematical Modeling,
IRL 2807 CNRS-UCHILE,
Facultad de Ciencias F\'isicas y Matem\'aticas,
Universidad de Chile,
Santiago, Chile.}}   
\date{25 January, 2023}
\maketitle

\begin{abstract}
We study a model for the spread of an infectious disease which
incorporates spatial and temporal effects. 
The model is a delayed multi-type branching process in
which types represent geographic regions while infected individuals 
reproduce offspring during a finite time
interval and have convalescence times and random death/recovery
outcomes.
We give simple expressions for the 
limit of the geometrically weighted mean evolution of
the process.
\end{abstract}

\bigskip

{\small\noindent

}

\bigskip

\section{Introduction} 

\noindent Inspired by agent-based simulation
studies of the spread of SARS-CoV2 such as that reported in
\cite{ferguson&etal2020}, our aim is to study a model for the
epidemic spread of an infectious agent which possesses 
levels of infectiousness that
vary both spatially and according to the time 
elapsed since infection within a fixed finite
time window and has a finite, but not necessarily bounded, convalescence time. 

\medskip

\noindent Our model is based on discrete-time branching processes (we
shall write ${\bp}$ for branching process in future), which find a
natural application in describing disease spread, for instance,
see~\cite{becker, Yanev&etal2020}.
More precisely, the model is a cross between a multi-type $\bp$ and a delayed
$\bp$ sporting the following characteristics. Individuals who contract an
 infectious disease are only contagious
for a short time immediately following their infection, say, $D$ days.
Thus, an individual born (infected) at time~$s$ has the opportunity to
reproduce offspring (cause infections) at times $\{s+d: 1\le d \le D\}$, and the number
of offspring born at time $s+d$ follows a law that depends on the time
offset $d$.
Individuals may continue suffering the effects of the disease long after
they have ceased to be infectious and the model assigns each individual
in the population a random lifetime (convalescence time). This time,
which is not necessarily bounded, may be zero, signifying that the
individual is asymptomatic. Symptomatic Individuals either die or
recover according to a Bernoulli random variable at the end of their
lifetimes. If they die before~$D$ days have lapsed, they cease to be
infectious. Asymptomatic individuals (with zero lifetimes) do not die,
but remain contagious for~$D$ days following infection.
Spatial disparity is captured by including a multi-type component in
which types represent regions and the numbers of offspring of all types
produced by individuals of distinct types can have different
distributions.
The model plays host to three closely related processes that evolve in
time: the offspring process $\X$, 
the asymptomatic population size process~$\Y$ and the
symptomatic population size process~$\Z$. 
In epidemiological terms, $\X$ captures disease incidence while $\Y+\Z$
models the prevalence.  

\medskip

\noindent The class of models we study here is much more restrictive
than the C-M-J $\bp$'s in which individuals reproduce according to a
random point process and are alive during random time intervals. These
are described and comprehensively studied in a large body of work (see
\cite{Crump&Mode1968,Crump&Mode1969,crump1970,
Doney1972,Doney1976,nerman1981}) that encompasses more general
frameworks than that considered here and they mostly deal with the
continuous-time setting. This literature examines the limiting behavior
of processes weighted by an exponential function of the Malthusian 
parameter in great detail and
has the establishment of long-term mean behavior, convergence in
distribution (see \cite[Theorem~2]{Doney1976}) and conditions for
a.s.~convergence (see \cite[Theorem~2]{Doney1976} and \cite[Proposition
~1.1]{nerman1981}) as notable achievements. In particular, concerning
the mean limits of processes, there are analytical procedures for
obtaining these (see \cite[Lemma~2 and Proposition~2]{Doney1976} and
\cite[Theorem~4.1]{crump1970}) which are general, but difficult to
employ in practical modeling applications in epidemiology where explicit
estimates are desirable.

\medskip

\noindent The paper is arranged as follows. Section~\ref{sec:delayed}
defines delayed multi-type ${\bp}$'s and describes the offspring
distribution specific to the model studied here. It also derives the
evolution equations for the mean population sizes of the~$\X$, $\Y$
and~$\Z$ processes, which are governed by a family of mean matrices at
different delay times. The $(i,j)$ element of the $d^{\rm th}$ mean
matrix indicates the average number of $j$-type individuals an $i$-type
individual infects a time~$d$ after contracting the disease.
Section \ref{sec: main} then presents the main result, 
Proposition~\ref{proplimits}, which gives simple-to-compute analytical
expressions for the limits of the mean population sizes of the three
processes weighted according to an exponential of the Malthusian parameter. 
These are analogous to results for classical multi-type ${\bp}$'s.

\medskip

\section{Delayed multi-type branching processes}
\label{sec:delayed}

\noindent In this model, each individual is born at some
time $s\in \NN_0=\{0,1,2,..\}$ and generates offspring 
independently of all other
individuals. The offspring can be of any type $i\in I$
and each is born within a finite set of times
$\aD\subset\NN=\{1,2,\ldots\}$ following the birth of the parent. Thus,
an individual born at $s$ generates offspring at times in $s+\aD$.
The number of offspring of each type in~$I$ produced by the same
parent at different ages in~$\aD$ are independent.
If $\min \aD>1$ there
is a latency period during which an individual is not yet contagious.
One can assume that g.c.d.$\,\aD\, =1$ and we write $D=\max \, \aD$ for
the maximum delay. Our interest is in the case $|\aD|>1$.

\medskip

\subsection{The random structure}

\noindent The set of nodes ${\aI}$ represents
all the potential individuals involved in the process.
A node $b\in {\aI}$ is identified by $b=(a;i,t,l)$
where $a$ is its parent node, $i\in I$ is its type, $t$ is its
time of birth and $l$ enumerates the nodes born to
parent $a$ of type $i$ at time~$t$. When $t=0$, $b$ is a root
and we set $a=\emptyset$. Each node $b$ gives rise to   
a set of nodes $(b;j,s,h)$ for $j\in I$, $s=t+d$
for $d\in \aD$ and $h\in \NN$.
Let~${\aI}^i$ denote the set of all type~$i$ nodes.

\medskip

\noindent Next, we associate with $b\in {\aI}$ independent random
elements $\bz(b)$ and $({\eL},\varepsilon)(b)$ which are also mutually 
independent.
Variables $\bz(b)$  
of the same type~$i$ are identically distributed for
$b\in {\aI}^i$.
For each $b=(a;i,t,l)$, $\bz(b)=(\bz_{a;t,l,d}^{i,j}: d\in \aD, j\in I)$
is a vector of random variables where
${\bz}_{a;t,l,d}^{i,j}$ is the number of potential offspring of type $j$
born to~$b$ at time $t+d$ and its distribution only depends on~$d$ and
$(i,j)$. We shall let ${\bz}_d^{i,j}$ denote a random variable having
this distribution.
The variables $({\eL},\varepsilon)(b)$ are identically distributed for
all $b\in {\aI}$. Variable ${\eL}(b)$ takes values in $\NN_0$ 
and indicates the lifetime of~$b$ while $\varepsilon(b)$ may depend on
${\eL}(b)$, taking the value~$1$ or~$0$ according as~$b$ recovers or
dies respectively after time ${\eL}(b)$ has elapsed.
If ${\eL}(b)=0$, $b$ is asymptomatic and cannot die, so we set
$\varepsilon(b)=1$. In contrast, if $\eL(b)>0$, then $b$ is symptomatic
and able to produce offspring at times in~$\aD$ up to~$D$, except if it dies
before~$D$. Here, we use `die' to mean a shift to a non-reproductive
state like actual death or isolation, as in \cite{Yanev&etal2020}.
Note that symptomatic individuals that recover before time~$D$
remain contagious and able to produce offspring at all times in~$\aD$ up 
to~$D$.
Let $({\eL},\varepsilon)$ denote a random element with the same law as
$({\eL},\varepsilon)(b)$. To avoid the possibility of having
trivial dynamics, we assume
$\PP(\eL<\infty)=1$ and $\PP(\eL>0)>0$,
and when modelling asymptomatic individuals we assume
$\PP(\eL=0)>0$.

\medskip

\noindent We now fix a random realization on ${\aI}^{{\NN_0}}$.     
In order for an individual $b=(a;t,l,d)$ to produce offspring as
described above, one requires that it has not yet died. So, the total
number of offspring of type $j$ born at time $t+d$ will be
$$
\xi_{a;t,l,d}^{i,j}={\bz}_{a;t,l,d}^{i,j}
\left(1-{\ind}\bigl(\eL(a;i,t,l)\le
d,\varepsilon(a;i,t,l)=0\bigr)\right).
$$
This ensures no offspring 
are produced after the death of the parent and fixes a dependence on the 
pair $(\eL,\varepsilon)$ specific to the individual.
The distribution of $\xi_{a;t,l,d}^{i,j}$
only depends on $d$ and $(i,j)$, and we will use
$\xi_d^{i,j}$ to denote a random variable having this distribution.
The individuals so generated, also denoted by $b\in {\aI}$ in what follows,
are identified by a triplet $(i,s,l)$ where $i$ is the type, $s$  
the time of birth and $l$ enumerates all
the individuals of type $i$ born at time $s$ and, from now on, we
dispense with the parent and simply write
${\bz}_{s,l,d}^{i,j}$ and $\xi_{s,l,d}^{i,j}$.
When an individual $b=(i,s,l)$ is generated by the process,
then it is ill and manifests symptoms in the time interval $[s,s+\eL-1]$
when $\eL>0$. In
contrast,if ${\eL}=0$, then once infected the individual is
asymptomatic during the time interval $[s,s+D]$
but is not counted as being ill. 
This means that asymptomatic individuals are able to infect others and
produce offspring during the interval $[s,s+D]$ while symptomatic
individuals can only do so during $[s,s+\min(\eL-1,D)]$.

\medskip

\noindent We should mention that having offspring born at a fixed set of
delays covers the case where the birth times of an individual's
offspring  are random and bounded. In fact, delays in a  bounded random time span~$\A$ can be viewed within the deterministic
framework as follows:
replace the number of $j$-type offspring produced by a $i$-type
individual at delay~$d$ by $\xi_d^{i,j}\ones(\A\ge d)$.

\medskip

\subsection{The processes}

\noindent The offspring process
${\X}(s)=\left({\X}_j(s): j\in I\right)$ is defined by
${\X}(s)=\zeros$ for $s<0$ and
\begin{equation}
\label{eqXmo1}
{\X}_j(s)=\ind(j=i_0, s=0)+\sum_{i\in I} \sum_{d\in \aD}
\sum_{l=1}^{{\X}_i(s-d)}\xi_{s-d,l,d}^{i,j} \text{ for } s\in\NN_0.
\end{equation}
The initial condition is a single type~$i_0$ individual
(which can be written as ${\X}(0)={\fe}_{i_0}$) and ${\X}_j(s)$ is
the number of type~$j$ offspring born at time~$s$ for $s>0$.
Since ${\X}=({\X}(s): s\ge 0)$ only counts offspring,
there is the implicit assumption that individuals
live for a single unit of time.

\medskip

\noindent Each individual $b=(i,t,l)$ has a lifetime ${\eL}(i,t,l)$
distributed as ${\eL}$ during which it is considered to be ill.
Define ${\Z}(s)=\left({\Z}_j(s): j\in I\right)$ to be the number of
ill (symptomatic) individuals of each type at time~$s$. Now, recall that
the set $\{(j, s,l):
l=1,\ldots,{\X}_j(s)\}$ enumerates the type~$j$ offspring born at
time~$s\ge 1$. For $s=0$, $(i_0,0,1)$ denotes the
initial individual. We set ${\Z}(s)=\zeros$ for $s<0$. 
Since an individual with ${\eL}(i,t,l)=0$ is never ill one has
\begin{equation}
\label{eqXmo2}
{\Z}_j(s)= \sum_{c=0}^s
\sum_{l=1}^{{\X}_j(s-c)} \ind({\eL}(j,s-c,l)>c),\; s\in \NN_0.
\end{equation}
We shall call ${\Z}=\left({\Z}(s):
s\in\NN_0\right)$ a delayed multi-type ${\bp}$.

\medskip

\noindent Individuals for whom $\eL=0$ are not counted as ill by the
model. They are asymptomatic for~$D$ time units and are
able to infect others during that time.
The process of asymptomatic cases present at each time is
${\Y}(s)=(\Y_j(s): j\in I)$ where
${\Y}(s)=\zeros$ for $s<0$ and
\begin{equation}
\label{mediaasi}
\Y_j(s) = \sum_{c=0}^D \sum_{l=1}^{\X_j(s-c)} \ind(\eL(j,s-c,l)=0),
\; s\in \NN_0.
\end{equation}
In the epidemiological setting, $\X$ models disease incidence while $\Y+\Z$
gives the prevalence.  

\medskip

\noindent The three processes $\Z$, $\X$ and ${\Y}$ become extinct
together almost surely or they all have some positive probability of not
dying out. In other words, the extinction time of~$\X$,
$T^{\X}=\inf\{t\ge 0: \sum_{j\in I}{\X}_j(t+s)=0 \text{ for all } 
s\ge0\}$,
and the analogously defined $T^{\Z}$ and $T^{\Y}$ satisfy:
$$
\PP(T^{\Z}<\infty)=1 \;\Longleftrightarrow\; \PP(T^{\X}<\infty)=1   
\;\Longleftrightarrow\; \PP(T^{\Y}<\infty)=1.
$$
The first equivalence is a direct consequence of the fact
that all individuals $b$ generated during this process have an
almost surely finite lifetime $\eL(b)$.
For the second equivalence, the implication ($\Rightarrow$)
follows from $T^\Y\le T^\X+D<\infty$ a.s. and conversely,
if $\Y$ becomes extinct and ${\X}$ is not extinct, one gets
$\PP\left(T^{\Y}<\infty, \exists t_n \to \infty,
\exists b_n=(i_n,t_n,1) \right)>0$.
Since lifetimes are identical and independent of
all other variables, imposing the condition that $\eL(b_n)=0$ merely
thins the set of individuals born after time $T^\Y$, so $\X$ must become
extinct.

\medskip

\noindent When averaging the processes $\Z$ and~$\X$ with 
respect to $(\varepsilon(b): b\in \aI)$, one obtains particular classes
of the processes considered in \cite{Doney1976} and \cite{nerman1981}. 
The individuals of 
${\X}$ have unit lifetimes, reproduce at time offsets in $\aD$ and
$i$-type individuals produce $j$-type offspring using a copy of 
$\xi_d^{i,j}$.
Similarly, individuals of ${\Z}$ have lifetimes given by
${\eL}>0$, reproduce at times in $\{d\in \aD: d<{\eL} \}$ and $i$-type
individuals reproduce $j$-type offspring using a copy of $\xi_d^{i,j}$.

\medskip

\noindent Let us average the law of $\xi_{t,l,d}^{i,j}$
over $(\eL,\varepsilon)$.
This gives us the offspring law of a contagious individual
when there is no information available about lifetime or
recovery $(\eL, \varepsilon)$.
The law $\bigl(p_d^{i,j}(n): n\in\NN_0\bigr)$ only depends on $i$,
$j$ and~$d$. Due to the   
independence between $\bz(b)$ and $(\eL,\varepsilon)(b)$,
$p_d^{i,j}(n)=\EE_{(\eL,\varepsilon)}\left(\PP\bigl(\xi_d^{i,j}=n
\bigr)\right)$, $n\ge 0$, takes the form
\begin{equation*}
p_d^{i,j}(n)=\begin{cases}
\PP({\bz}_d^{i,j}=n) \bigl(1-\PP({\eL}\le d,\varepsilon=0)\bigr),
& \text{ if } n>0, \\
\PP({\bz}_d^{i,j}=0) \bigl(1-\PP(\eL\le d,\varepsilon=0)\bigr)
+ \PP(\eL\le d,\varepsilon=0), & \text{ if } n=0.
\end{cases} 
\end{equation*}
Naturally, $\sum_{n\ge 0}p_d^{i,j}(n)=1$. The mean number of offspring of
type~$j$ produced by an individual of type~$i$ and age~$d$ is then given by
$$
M_d(i,j)=\EE\left(\xi_d^{i,j}\right)
=\sum_{n\ge 0} n \, p_d^{i,j}(n) ,
\text{ for } i,j\in I, d\in \aD.
$$
Let $M_d=\left(M_d(i,j): i,j\in I\right)$. We assume $M_d$ is irreducible
for all $d\in \aD$. By
convention, set $M_d=\zeros$ and $p_d^{i,j}(n)=\ind(n=0)$ 
for any $d\not\in \aD$. 

\medskip

\noindent  We will denote the expected value when starting with a single
individual of type $i_0$ by $\EE_{i_0}$, but if there is no confusion we
shall simply write $\EE$.
Set $\EE({\X}(s))=\zeros$ for $s<0$. Since ${\X}(s-d))$ is independent
of the sequence of variables $(\xi_{s-d,l,d}^{i,j})$ by Wald's
equation~\cite{Wald1944}, taking expectations on both sides of
~(\ref{eqXmo1}) gives
\begin{equation}
\label{evolutionmod1}
\EE_{i_0}({\X}(s))' = {\fe}_{i_0}'\ind(s=0)
+\sum_{d\in \aD} \EE_{i_0}({\X}(s-d))' M_d, \; s\in \NN_0.
\end{equation}
To evaluate the expected number of individuals of each type ill at
time~$s$, $\EE_{i_0}({\Z}_j(s))$, we can use $M_d=\zeros$ for $d\not\in
\aD$ together with relation~(\ref{eqXmo2}) and Wald's
Equation~\cite{Wald1944} once again to obtain
\begin{equation}
\label{evolutionmod2}
\EE_{i_0}({\Z}_j(s))
= \sum_{c=0}^s \EE_{i_0}\left( \sum_{l=1}^{{\X}_j(s-c)}
\ind\bigl(\eL(j,s-c,l)> c\bigr) \right)
= \sum_{c=0}^s \EE_{i_0}(\X_j(s-c)) \PP(\eL> c).
\end{equation}
Now, from (\ref{mediaasi}) and the independence properties of ${\eL}$ and ${\X},$ the 
expected number of asymptomatic cases is given by
\begin{equation}
\label{evolutionmod3}
\EE_{i_0}(\Y(s))' \!=\! \PP(\eL=0) \sum_{d=0}^D \EE_{i_0}(\X(s-d))', \;s\! \in\! \NN_0.
\end{equation}
  
\section{The encoding and the main result}
\label{sec: main}

\noindent Following Definition~2 in~\cite{Doney1976} where the
Malthusian  parameter is given for general delayed multi-type processes,
we have that the 
exponential of the Malthusian parameter of~${\X}$,
say~$\lambda$, is uniquely defined by
\begin{equation}
\label{zzurho1}
\lambda \hbox{ is the unique value for which 
the Perron-Frobenius eigenvalue of } \sum_{d\in \aD}\lambda^{-d} M_d \text{ is } 1.
\end{equation}

\medskip

\noindent We are interested in $\lim\limits_{s\to \infty} 
\lambda^{-s}\EE(({\X}(s))$. 
First, notice that ${\X}$ is equally distributed as a delayed ${\bp}$
${\hX}$ whose definition includes no concept of lifetime and where 
each individual $(i,t,l)$ generates $\hxi_{t,l;d}^{i,j}$
offspring of type~$j$ at time $t+d$ according to law
$p_d^{i,j}$. So, $\hX$ satisfies 
${\hX}(s)=\zeros$ for $s<0$ and
${\hX}_j(s)=\ind(j=i_0, s=0)+\sum_{i\in I} \sum_{d\in \aD}
\sum_{l=1}^{{\hX}_i(s-d)}\hxi_{s-d,l,d}^{i,j}$ for $s\in\NN_0$.
The mean number of offspring in ${\hX}$ and ${\X}$ is given by
$M_d^{i,j}=\sum_{n\ge 0} n p^{i,j}_d(n)$, and
$\EE({\hX}(s))=\EE({\X}(s)), \, s\in \NN_0$. Now we will encode
${\hX}$ by a multi-type ${\bp}$.

\medskip

\noindent Set $[D]=\{1,\ldots,D\}$ and consider the new set of types 
${\hI}=[D]\times I$. Then, ${\hX}$ may be viewed as the following 
multi-type ${\bp}$ ${\hZ}$ on the set of types ${\hI}$ defined by:
\begin{equation*}
{\hZ}(s)=({\hZ}_{d,j}(s): (d,j)\in {\hI}), \quad s\ge 0,
\text{ with } 
{\hZ}_{d,j}(s)={\hX}_j(s+1-d).
\end{equation*}
So, ${\hZ}_{d,j}(0)=\ind(d=1)\ind(j=i_0)$.
For $1< d\le D$ one has
${\hZ}_{d,j}(s+1)={\hZ}_{d-1,j}(s)$
while for $d=1$ one obtains from (\ref{eqXmo1}) that
$$
{\hZ}_{1,j}(s+1))={\hX}_j(s+1)=\sum_{i\in I} \sum_{e\in \aD}
\sum_{l=1}^{{\hX}_i(s+1-e)} {\hxi}_{s+1-e,l,e}^{i,j}=
\sum_{i\in I} \sum_{e\in \aD}
\sum_{l=1}^{{\hZ}_{e,i}(s)}{\hxi}_{s+1-e,l,e}^{i,j}.
$$
The mean matrix of~$\hZ$ is given by 
${\tM}((e,i),(d,j))=\ind(i=j)\ind(e=d-1)$
if $d>1$. When $d=1$, one has ${\tM}((e,i),(1,j))=0$ if $e\notin \aD$ and
${\tM}((e,i),(1,j))=M_e(i,j)$ otherwise.
One can check that the matrix ${\tM}$ is irreducible.
Let $\hnu=(\hnu(e,i):(e,i)\in \hI)$ be a left eigenvector
of~$\tM$ corresponding to its Perron-Frobenius eigenvalue ${\rho}$ and let
$\hnu_e=(\hnu(e,i): i\in I)$ for $e\in [D]$. For $d>1$, one has
${\rho}{\hnu}_d={\hnu}_{d-1}$ and so
$$
{\hnu}_d={\rho}^{-(d-1)} {\hnu}_1,\, d\in [D]. 
$$
Hence,
$$
{\rho}{\hnu}_1(j)
=\sum_{e\in \aD}\sum_{i\in I} {\rho}^{-(e-1)} {\hnu}_1(i)M_e(i,j)
=\sum_{i\in I} {\hnu}_1(i)
\left(\sum_{e\in \aD}{\rho}^{-(e-1)} M_e(i,j)\right),
$$
from which it follows that ${\hnu}_1$ satisfies $\hnu'_1 M_\rho=\hnu'_1$, where
$$
M_\rho=\sum_{e\in \aD}{\rho}^{-e} M_e.
$$
Therefore$1$ is the Perron-Frobenius eigenvalue of $M_\rho$ and
it follows that~$\lambda$ defined in (\ref{zzurho1}) is also the
Perron-Frobenius eigenvalue of $\tM$. 
Similarly take $\tth=(\tth(e,i):(e,i)\in \hI)$ to be a right Perron-Frobenius eigenvector
of~$\tM$ so that $\tM {\tth}=\lambda \tth$. Write 
${\tth}_e=(\tth(e,j): j\in  I)$ for $e\in [D]$. Then
$M_e {\tth}_1=\lambda {\tth}_e-{\tth}_{e+1} \ind(e<D)$ and 
hence $\lambda^{-e} M_e {\tth}_1=\lambda^{-(e-1)} {\tth}_e-\lambda^{-e}{\tth}_{e+1} 
\ind(e<D)$. By summing over $e=1,\ldots,D$, we obtain $M_\lambda 
{\tth}_1={\tth}_1$,
so that ${\tth}_1$ is the eigenvector corresponding to eigenvalue $1$ of 
the matrix $M_\lambda$. Also we can iterate equality $M_e {\tth}_1=\lambda 
{\tth}_e-{\tth}_{e+1} \ind(e<D)$ to get 
$$
{\tth}_d=\sum_{e=d}^D {\lambda}^{e-d+1}M_e {\tth}_1 \text{ for } d\in [D].
$$
(Recall that $M_d=\zeros$ when $d\not\in \aD$).

\medskip

\noindent Since $({\hZ}(s): s\ge 0)$ is a multi-type ${\bp}$  
with mean matrix $\tM$, we have 
$$
\EE(({\hZ}(s))=\EE({\hZ}(0))\tM^s. 
$$
Also ${\hZ}_{d,j}(0)=\ind(d=1)\ind(j=i_0)$ which means
$\EE({\hZ}(0))={\fe}_{1,i_0}$. 

\medskip

\noindent It is well known that if $A$ is an irreducible non-negative matrix
with Perron-Frobenius eigenvalue $\rho_A$ and corresponding left and right  eigenvectors
normalized so that $\nu' h=1$, Then, for any norm $\norm{\cdot}$,
there exists $C<\infty$ and $\delta\in (0,1)$
such that $\norm{\rho^{-s} A^s-h\nu'}\le C \delta^s$ for all $s\ge 0$ and
$\lim\limits_{s\to \infty} \rho^{-s} A^s= h \nu'$ componentwise
(See Chapter~1 of~\cite{es81}).

\medskip

\noindent So, when we normalize $\hnu$ and ${\tth}$ such that $\hnu' {\tth}=1$ 
we have $\lim\limits_{s\to \infty} \lambda^{-s}  {\tM}^s= {\tth} {\hnu}'$ 
componentwise. Therefore,
$$
\lim\limits_{s\to \infty} \lambda^{-s}\EE(({\hZ}(s))
= {\fe}_{1,i_0}{\tth} {\hnu}'. 
$$
Using $\EE(\X(s)) = 
\EE(\hX(s)) =\EE(\hZ_{1,\cdot}(s))$, we get
\begin{equation}
\label{limnorm2}
\lim\limits_{s\to \infty} \lambda^{-s}\EE(({\X}(s))= {\fe}_{i_0}{\tth}_1 {\hnu}'_1.
\end{equation}

\noindent Next we compute the condition $\hnu'{\tth}=1$ in terms
of $\hnu_1$ and ${\tth}_1$. We have
\begin{eqnarray}
\nonumber
1 &=& \hnu'{\tth}
= \sum_{d=1}^D \lambda^{-(d-1)}\hnu_1' 
(\sum_{e=d}^D \lambda^{-(1+e-d)}M_e){\tth}_1
=\hnu_1' \sum_{d=1}^D (\sum_{e=d}^D \lambda^{-e}M_e){\tth}_1\\
\label{normalix}
&=& \hnu_1' (\sum_{e=1}^D\sum_{d=1}^e \lambda^{-e}M_e){\tth}_1
= \hnu_1' \left(\sum_{d=1}^D d\, \lambda^{-d}\, M_d\right){\tth}_1.
\end{eqnarray}  

\medskip

\noindent Let $\nu$ and $h$ be 
the left and right eigenvectors 
$\hnu_1$ and ${\tth}_1$ of $M_\lambda$ but normalized so that
$\nu' h=1$. Define
\begin{equation*}
\mu= \nu' \left(\sum_{d=1}^D d\, \lambda^{-d} \, M_d\right) h.
\end{equation*}

\begin{proposition}
\label{proplimits}
The limit mean for the offspring process is
\begin{equation}
\label{limnorm3}
\lim\limits_{s\to \infty} \lambda^{-s}\EE(({\X}(s))
= \mu^{-1}{\fe}_{i_0} h \nu'
\end{equation}
while the limit mean for the asymptomatic individuals is
\begin{equation}
\label{limnorm4}
\lim_{s\to\infty} \lambda^{-s}\EE_{i_0}(\Y(s))' 
= \mu^{-1} \PP(\eL=0) \left( \sum_{d=0}^D \lambda^{d} \right)
{\fe}_{i_0} h \nu'.
\end{equation}
Some additional assumptions are required to handle process~$\Z$: If the
processes are subcritical 
($\lambda<1$), assume that
$\EE\left(\lambda^{-\eL}\right)<\infty$, while if they are critical
($\lambda=1$), assume that $\EE(\eL)<\infty$.
Then, the long-term behavior of the mean number of symptomatic
individuals is given by
\begin{equation*}
\lim_{s\to\infty} \lambda^{-s}\EE_{i_0}({\Z}_j(s))
=\mu^{-1} \left ( \sum_{c=0}^\infty \lambda^{-c}\PP(\eL> c). \right) 
{\fe}_{i_0} h \nu'.
\end{equation*}
\end{proposition}

\begin{proof}
Equation (\ref{limnorm3}) follows directly from (\ref{limnorm2}) and
(\ref{normalix}). This equality together with (\ref{evolutionmod3}) then
leads to (\ref{limnorm4}). Next consider the case for symptomatic
individuals. From (\ref{evolutionmod2}), we have
$$
\EE_{i_0}({\Z}_j(s)) = \sum_{c=0}^s \EE_{i_0}(\X_j(s-c)) \PP(\eL> c).
$$
Multiplying both sides of this by $\lambda^{-s}$ and taking limits as $s\to\infty$ yields
\begin{eqnarray*}
&{}&\lim_{s\to\infty} \lambda^{-s}\EE_{i_0}({\Z}_j(s))
= \lim_{s\to\infty} \lambda^{-s}\sum_{c=0}^s
\EE_{i_0}(\X_j(s-c)) \PP(\eL> c)\\
&{}& = \lim_{s\to\infty}
\sum_{c=0}^s \left(\lambda^{-(s-c)}\EE_{i_0}(\X_j(s-c)) \right)
\lambda^{-c}\PP(\eL> c)\\
&{}& = \lim_{s\to\infty}
\sum_{c=0}^\infty \left(\lambda^{-(s-c)}
\EE_{i_0}(\X_j(s-c))\ind(s\ge c) \right) \lambda^{-c}\PP(\eL> c).
\end{eqnarray*}
Next,
$$
\lambda^{-(s-c)}\EE_{i_0}(\X_j(s-c))\ind(s\ge c) =
\lambda^{-(s-c)}\EE_{i_0}(\X_j(s-c))
$$
because $\X_j(s)=0$ for $s<0$. We have already seen that
$$
\lim_{s\to\infty} \lambda^{-(s-c)}\EE_{i_0}(\X_j(s-c)) 
= \mu^{-1}{\fe}_{i_0}h
\nu',
$$
for $c\ge0$. Since
$\lim\limits_{s\to\infty}\lambda^{-(s-c)}\EE_{i_0}(\X_j(s-c))$ is
finite,  we have
$$
\sup\{\lambda^{-s}\EE_{i_0}(\X_j(s)): s\ge0\} < \infty,
$$
so there is some $K>0$ that bounds
$\lambda^{-(s-c)}\EE_{i_0}(\X_j(s-c))$ for all $s\ge 0, 0\le c\le s$.
If~$\Z$ is supercritical, then $\lambda>1$ and
$$
\sum_{c=0}^\infty \left(\lambda^{-(s-c)}\EE_{i_0}(\X_j(s-c)) \right)
\lambda^{-c}\PP(\eL> c)
\le K \sum_{c=0}^\infty \lambda^{-c} \PP(\eL> c)  < \infty.
$$
However, this also holds under the additional assumptions made for the critical
and subcritical cases. Consequently, the dominated convergence theorem
can be used to justify exchanging the limit with the summation and we
obtain
\begin{align*}
& \lim_{s\to\infty} \lambda^{-s}\EE_{i_0}({\Z}_j(s))
= \lim_{s\to\infty} \sum_{c=0}^s \left(\lambda^{-(s-c)}
\EE_{i_0}(\X_j(s-c)) \right) \lambda^{-c}\PP(\eL> c)\\
&= \sum_{c=0}^\infty \left(\lim_{s\to\infty} \lambda^{-(s-c)}
\EE_{i_0}(\X_j(s-c)) \right) \lambda^{-c}\PP(\eL> c)\\
& = \mu^{-1} \left( \sum_{c=0}^\infty \PP(\eL> c) \lambda^{-c} \right) 
{\fe}_{i_0} h \nu'.
\tag*{\qedhere}
\end{align*}
\end{proof}

\medskip

\begin{corollary}
The normalized left eigenvector $\nu'$ 
is the limit of types for the
processes~$\X$, $\Y$ and $\Z$, that is,
\begin{equation*}
\lim_{s\to\infty} \frac{\EE(\X(s))'}{\EE(\X(s))'\ones}
=\lim_{s\to\infty} \frac{\EE(\Y(s))'}{\EE(\Y(s))'\ones}
=\lim_{s\to\infty} \frac{\EE(\Z(s))'}{\EE(\Z(s))'\ones}
= \nu'.
\end{equation*}
\end{corollary}

\begin{proof}
This follows straightforwardly from Proposition \ref{proplimits}
since $\nu' \ones=1$.
As the calculation of the limits is the same for all three processes, 
we will only present those for~$\Z$.
\begin{equation*}
 \lim_{s\to\infty} \frac{\EE(\Z(s))'}{\EE(\Z(s))'\ones}
=\frac{\mu^{-1}(\EE({\Z}(0))' h)
\left(\sum_{c\ge 0} \PP(\eL> c) \lambda^{-c}\right) \nu'}
{\mu^{-1}(\EE({\Z}(0))' h)
\left(\sum_{c\ge 0} \PP(\eL> c) \lambda^{-c}\right) \nu'\ones}
= \nu' \,.
\tag*{\qedhere}
\end{equation*}
\end{proof}

\medskip

\begin{remark}
Analogous to what happens in the case of a multi-type ${\bp}$, $\nu'$ is
a stationary distribution for the mean evolution of types for the
process~$\X$. In fact, $\nu'$ is the left eigenvector corresponding to the eigenvalue~$1$
of $M_{\lambda}=\sum_{d\in\aD} {\lambda}^{-d}M_d$ normalized to sum to
unity. So when it is taken as the initial distribution one gets
$\EE({\X}(s))' =\nu' \lambda^s$ for $s=0,1,\ldots,D-1$
and (\ref{evolutionmod1}) yields
$\EE({{\X}}(s))' = \sum_{d\in \aD} \nu' M_d \lambda^{s-d}
= \lambda^s \nu'\sum_{d\in \aD} \lambda^{-d} M_d 
= \nu' \lambda^s$.  
Therefore, $\EE({\X}(s))')/(\EE({\X}(s))'\ones) = \nu'$
for all $s\ge 0$.
\end{remark}

\medskip

\section*{Acknowledgments}

This work was supported by the Center for Mathematical Modeling ANID Basal
Projects ACE210010 and FB210005.


\end{document}